\documentclass[11pt]{amsart}

\usepackage{amsmath}
\usepackage{amssymb}
\usepackage{amsfonts}
\newtheorem{theorem}{Theorem}
\theoremstyle{plain}

\newtheorem{corollary}{Corollary}

\newtheorem{definition}{Definition}

\newtheorem{lemma}{Lemma}

\newtheorem{proposition}{Proposition}

\numberwithin{equation}{section}
\def\<{{\langle}}
\def\>{{\rangle}}
\def\S{{\mathbb{S}}}

\begin{document}
\title[Exceptional minimal surfaces in spheres]
{Exceptional minimal surfaces in spheres}
\author{Theodoros Vlachos}
\address{Department of Mathematics, University of Ioannina, 45110 Ioannina,
Greece}
\email{tvlachos@uoi.gr}
\subjclass[2000]{Primary 53A10; Secondary 53C42}
\keywords{Minimal surface, sphere, Hopf differentials, exceptional
surface, pseudoholomorphic curves in $\S^{6}$, polar map}

\maketitle

\begin{abstract}
We  study a class  of exceptional minimal surfaces in spheres for which all  Hopf differentials are holomorphic. 
Extending results of Eschenburg and Tribuzy \cite{ET0}, we obtain a description of  exceptional surfaces  in terms of a set of
absolute value type functions, the $a$-invariants, that determine
the geometry of the higher order curvature ellipses and 
satisfy certain Ricci-type conditions. We show that the $a$-invariants determine these surfaces 
 up to a  multiparameter family of isometric minimal deformations, where the number of the parameters
is precisely the number of non-vanishing Hopf differentials.
We give applications to superconformal surfaces and pseudoholomorphic curves 
in the nearly K\"{a}hler sphere $\S^{6}$. Moreover, we study superconformal surfaces in odd dimensional spheres that are isometric to their polar and show a relation to pseudoholomorphic curves in $\S^{6}$.

\end{abstract}

\footnotetext[1]{The author was supported by the Alexander von Humboldt Foundation.}

\section{Introduction}
Minimal surfaces in the Euclidean space $\mathbb{R}^{3}$ are locally constructed via the Weierstrass representation. More 
generally,
any minimal surface in $\mathbb{R}^{n}$ is locally the real part of the integral of an isotropic curve in $\mathbb{C}^{n}.$ 
The study of minimal surfaces in spheres was initiated by Calabi in his seminal paper \cite{Cal} 
and follows a  completely different 
route. 

In this paper, we consider \textit{%
exceptional}\footnotetext[2]{%
We note that this terminology was used by Johnson \cite{J}, but his surfaces
are superconformal (resp. superminimal) according as they lie  in
an odd-(resp. even-) dimensional sphere with substantial codimension.} surfaces, a class of minimal surfaces in spheres for which 
certain invariants, the so called Hopf differentials, are holomorphic. 
 This class of 
surfaces, includes the superminimal ones, which are the minimal surfaces with  vanishing
 Hopf differentials, minimal
two-spheres \cite{Ch},
superconformal surfaces \cite{BPW, M1}, as well as pseudoholomorphic curves in the nearly K\"{a}hler sphere $\S^6$ introduced in \cite{B}. 
Besides flat minimal surfaces (cf. \cite{K1,K2}), Lawson's  surfaces, i.e.,
minimal surfaces that decompose as a direct sum of the associated minimal
surfaces in $\S^{3}$, are indeed exceptional (see \cite{V2}). These surfaces are related to Lawson's conjecture \cite{L} which
asserts that the only non-flat minimal surfaces in spheres that are locally isometric to minimal surfaces in $\S^3$ are Lawson's 
surfaces.

The Hopf differentials play the role that Frenet curvatures play for curves, in the sense that two isometric minimal surfaces with 
the same Hopf differentials are congruent \cite{V}. 
Using a null basis for each higher complexified normal  bundle,  we split the Hopf 
differentials into a product of two factors.
 The modulus of each factor
defines scalar invariants which we call $a$-\textit{invariants} and determine the geometry of the higher 
curvature ellipses.

Our aim is to give a complete description of exceptional surfaces in spheres
in terms of the $a$-invariants, and to characterise their induced metrics. 
In fact, we give an existence and uniqueness theorem for exceptional surfaces in terms of the $a$-invariants  in the spirit of \cite{EGT,ET0}.
 We prove that for each exceptional surface, the $a$-invariants  are of absolute value type functions in the sense of
 \cite{EGT,ET} and satisfy certain conditions. Conversely,  we show  that each set of 
absolute value type functions that fulfil these conditions  determines an exceptional surface  up 
to a multiparameter family of isometric minimal deformations, where the number of the parameters is precisely the number of non-vanishing Hopf differentials.

Several applications are provided. 
At first, we investigate superconformal surfaces that are isometric to their polar. The Gauss map of a minimal surface $M$ in $\S^3$ defines another minimal,
 the polar  \cite{L0}   of $M$, which is conformal to $M$.  The polar is also  defined for any minimal surface lying 
in an odd-dimensional sphere (cf. \cite{DF}). It was proved by Miyaoka \cite{M1} 
that the polar of any superconformal surface is again a superconformal surface. Moreover, this construction is dual, in 
the sense that taking the polar a second time produces the original surface. 

As an application of the main result, we study superconformal surfaces that are isometric to their polar,  which we briefly  call 
 \textit{self-dual surfaces}. 
 In contrast to the case of $\S^3$, where the Clifford torus is the only self-dual surface, there are lots of self-dual surfaces 
in high codimension. Flat superconformal surfaces and  Lawson's surfaces in  $\S^{8k+7}$ that are superconformal  are self-dual. 

The case of self-dual surfaces in $\S^5$ is quite interesting. We show that a 
superconformal surface in $\S^5$ is self-dual if and only if it is congruent to a
 pseudoholomorphic curve of $\S^6$ lying  in a totally geodesic $\S^5$. More generally, we  characterise  all 
 self-dual surfaces. 

Moreover, we give applications to pseudoholomorphic curves in the nearly K\"{a}hler 
sphere $\S^6$. We provide extrinsic characterizations in terms 
of the $a$-invariants, among the class of superconformal surfaces and reprove the intrinsic characterization given by Hashimoto \cite{H}. Finally, we give another short proof of the Lawson's conjecture for exceptional surfaces in 
spheres of odd dimension  \cite{V2}.

The paper is organised as follows: In section 2, we fix the notation and
give some preliminaries. In section 3, we consider the splitting of the Hopf differentials and introduce the  $a$-invariants.
 In section 4,  we prove the main result of the paper, namely that the $a$-invariants determine
all exceptional surfaces up to a multiparameter family. In section 5, we deal with pseudoholomorphic curves 
in the nearly K\"{a}hler sphere $\S^6$. Section 6 is devoted to the Ricci condition.
In Section 7, we investigate  self-dual surfaces. Finally some global formulas and topological restrictions are obtained. 

\section{Preliminaries}
Let $f:(M,ds^{2})\rightarrow \S^{n}$ be a minimal surface. Curves on $M$ through a point $p\in M$ have their first
derivatives on the tangent plane $T_{p}M$, but higher order derivatives will
have components normal to $f$. The $r$-th osculating space $T_{p}^{r}f$ of $f$ at $p$ is spanned by the
derivatives of order up to $r$ is called the  of $f$
at $p$ (cf. \cite{Ch,Sp}).  The $r$-th normal space  $N_{p}^{r}f$ is the orthogonal complement of $T_{p}^{r}f$ in $T_{p}^{r+1}f$ and has dimension $\leq 2$. 

A point $p$ is called generic (cf. \cite{Ch,A}) if dim$N_{p}^{r}f=2$ for any $%
r $, unless $r=m$ and $n=2m+1$, where dim$N_{p}^{m}f=1$. If $f$ is substantial in $%
\S^{n}$ in the sense of \cite{DR}, then the set of generic points is open and dense (cf. \cite[p. 96]{O}). Hereafter,
 we always assume that the minimal surfaces under consideration are substantial, unless otherwise stated. At generic points,
we can consider the $r$-th normal bundle $N^{r}f$ of $f$, with fibers $%
N_{p}^{r}f$, which is a plane bundle,
unless $n=2m+1$, where $N^{m}f$ is a line bundle.

 The $(r+1)$-th fundamental form $%
B_{r}$ is the $(r+1)$-linear tensor from $T_{p}M$ into $N_{p}^{r}f$, defined
by%
\begin{equation*}
B_{r}(X_{1},...,X_{r+1})=\pi _{r}\big(\overline{\nabla }_{\overline{X}%
_{1}}...\overline{\nabla }_{\overline{X}_{r}}\overline{X}_{r+1}\big),
\end{equation*}%
where $\pi _{r}$ is the orthogonal projection onto $N_{p}^{r}f$, $\overline{%
\nabla }$ is the Levi-Civit\'{a} connection of $\S^{n}$, and $\overline{X}%
_{1},...,\overline{X}_{r+1}$ are local vector fields that extend ${X}%
_{1},...,{X}_{r+1}$. It is well known that $B_{r}$ is symmetric (cf. \cite[%
p. 240]{Sp}) and $N_{p}^{r}f$ is spanned by the image of $B_{r}$. Clearly, $%
B_{1}$ is nothing but the second fundamental form.

We use the moving frame method and adopt the following convention on the
range of indices (the symbol $i$ is reserved for $\sqrt{-1}),$ unless
otherwise stated: 
\begin{equation*}
1\leq j,k\leq 2,{\ }3\leq \alpha,\beta \leq n,{\ }1\leq A,B,C\leq n,{\ }%
1\leq r,s,t\leq \Big[\frac{n}{2}\Big].
\end{equation*}%
Let $\{e_{A}\}$ be a local orthonormal frame field on $\S^{n}$, and let $\{\omega
_{A}\}$ be the coframe dual to $\{e_{A}\}$. The structure equations of $\S^{n}$ are 
\begin{eqnarray}\label{s1}
d\omega _{A} &=&\sum_{B}\omega _{AB}\wedge \omega _{B}, 
\end{eqnarray}
\begin{eqnarray}\label{s2}
d\omega _{AB} &=&\sum_{C}\omega _{AC}\wedge \omega _{CB}-\omega _{A}\wedge
\omega _{B},
\end{eqnarray}%
where the connection form $\omega _{AB}$ is given by $\omega
_{AB}(X)=\langle \overline{\nabla }_{X}e_{A},e_{B}\rangle $ and $\langle
.,.\rangle $ is the Riemannian metric on $\S^{n}$. We choose the frame such
that, restricted to $M$, $e_{j}$ is tangent and consequently $e_{\alpha }$
is normal to the surface. Then we have $\omega _{\alpha }=0$. By (\ref{s1}) and
Cartan's Lemma, we obtain
\begin{equation*}
\omega _{j\alpha }=\sum_{k}h_{jk}^{\alpha }\omega _{k},{\ }h_{jk}^{\alpha
}=h_{kj}^{\alpha }.
\end{equation*}%
The assumption that $f$ is minimal is equivalent to $h_{11}^{\alpha
}+h_{22}^{\alpha }=0$. Restricting equations (\ref{s1}) and (\ref{s2}) to $M,$ we obtain
the Cartan structure equations of $f.$

Hereafter we set $m:=[(n-1)/2]$, and  choose the
normal frame $e_{\alpha }$ such that $(e_{2r+1},e_{2r+2})$ is a frame field of $N^{r}f$ for any $r\leq m$. When $n=2m+1,$ $%
e_{2m+1}$ spans the fibers of $N^{m}f$. Then it is easy to see that (cf. 
\cite[Lemma 69]{Sp}) 
\begin{equation}\label{om}
\omega _{2r-1,\alpha }=\omega _{2r,\alpha }=0{\ }{\ }{\ }{\hbox {if}}{\ }{\ }%
{\ }\alpha >2r+2{\ }{\ }{\ }{\hbox {or}}{\ }{\ }{\ }\alpha <2r-3.
\end{equation}%
The components of the higher fundamental forms are given by 
\begin{equation*}
h_{1}^{\alpha }:=\langle
B_{r}(e_{1},...,e_{1}),e_{\alpha }\rangle, {\ } h_{2}^{\alpha }:=\langle
B_{r}(e_{1},...,e_{1},e_2),e_{\alpha }\rangle,
\end{equation*}%
where $\alpha =2r+1$ or $2r+2$. We use complex
vectors, and we put 
\begin{equation*}
H_{\alpha }=h_{1}^{\alpha }+ih_{2}^{\alpha },\text{ }E=e_{1}-ie_{2}{\ \ }%
\text{and }\varphi =\omega _{1}+i\omega _{2}.
\end{equation*}%
Then we have  (cf. \cite[p. 30]{Ch}): 
\begin{equation}\label{hi}
H_{2r+1}\omega _{2r+1,\alpha }+H_{2r+2}\omega _{2r+2,\alpha }=H_{\alpha }%
\overline{\varphi }\text{ \ for }\alpha =2r+3,2r+4,
\end{equation}
$0\leq r\leq m-1$, where $H_{1}=1,H_{2}=i,$ and when $n=2m+1$ 
\begin{equation*}
H_{2m-1}\omega _{2m-1,2m+1}+H_{2m}\omega _{2m,2m+1}=H_{2m+1}\overline{%
\varphi }.
\end{equation*}
The induced metric is $ds^{2}=\varphi \overline{\varphi }.$ From (\ref{s1}), we
find 
\begin{equation}\label{d}
d\varphi =-i\omega _{12}\wedge \varphi
\end{equation}%
and the Gaussian curvature $K$ is given by 
\begin{equation}\label{K}
d\omega _{12}=-\frac{i}{2}K\varphi \wedge \overline{\varphi }.
\end{equation}%

The $r$-\textit{th curvature ellipse},  at any point $%
p$ in $M$, is given by 
\begin{equation*}
\mathcal{E}\sb r(p)=\big\{B_{r}(X,...,X)\in N_{p}^{r}f:X\in T_{p}M,|X|=1%
\big\}.
\end{equation*}%
It is known (cf. \cite{Chen}) that $\mathcal{E}\sb r(p)$ is indeed an
ellipse (possibly degenerated). The $r$-\textit{th normal curvature} $K_{r}^{\perp }$
is defined by%
\begin{equation}\label{n}
K_{r}^{\perp }=i\left( H_{2r+1}{\overline{H}_{2r+2}}-{\overline{H}_{2r+1}}%
H_{2r+2}\right).
\end{equation}%
We note that the sign of $K_{r}^{\perp }$ depends on the orientation of the bundle $N^{r}f$. 
It is obvious that  $K_{r}^{\perp }(p)=0$ if
and only if ${\hbox {dim}}N_{p}^{r}f\leq 1$. 
It is not hard to verify that
$$
|K_{r}^{\perp }|={\frac{2}{\pi }}{\hbox {Area}}(\mathcal{E}\sb r),
$$
 or equivalently
\begin{equation}\label{n1}
|K_{r}^{\perp }|=2\kappa _{r}\mu _{r},
\end{equation}%
where   $\kappa _{r}\geq \mu
_{r}\geq 0$ are the length of the semi-axes of $\mathcal{E}\sb r$. 
For the sake of convenience, we also set $K_{0}^{\perp }=2.$ The length of $
B_{r}$ is given by%
\begin{equation}\label{l}
\left\Vert B_{r}\right\Vert ^{2}=2^{r}\big(|{H_{2r+1}}|^{2}+|{H_{2r+2}}|^{2}%
\big),
\end{equation}%
or equivalently (cf. \cite{A})%
\begin{equation}\label{l1}
\left\Vert B_{r}\right\Vert ^{2}=2^{r}(\kappa _{r}^{2}+\mu _{r}^{2}).
\end{equation}%

At points where $K_{r}^{\perp } \neq 0,r\geq 0,$  by using (\ref{n}) and (\ref{hi}), we find
\begin{equation}\label{01}
\omega _{2r+1,2r+3}=\frac{i}{K_{r}^{\perp }}\left( H_{2r+3}\overline{H}%
_{2r+2}\overline{\varphi }-H_{2r+2}\overline{H}_{2r+3}\varphi \right),
\end{equation}%
\begin{equation}\label{02}
\omega _{2r+2,2r+3}=\frac{i}{K_{r}^{\perp }}\left( H_{2r+1}\overline{H}%
_{2r+3}\varphi -H_{2r+3}\overline{H}_{2r+1}\overline{\varphi }\right),
\end{equation}%
\begin{equation}\label{03}
\omega _{2r+1,2r+4}=\frac{i}{K_{r}^{\perp }}\left( H_{2r+4}\overline{H}%
_{2r+2}\overline{\varphi }-H_{2r+2}\overline{H}_{2r+4}\varphi \right),
\end{equation}%
\begin{equation}\label{04}
\omega _{2r+2,2r+4}=\frac{i}{K_{r}^{\perp }}\left( H_{2r+1}\overline{H}%
_{2r+4}\varphi -H_{2r+4}\overline{H}_{2r+1}\overline{\varphi }\right).
\end{equation}%

Taking the exterior derivative of (\ref{hi}), and using (\ref{s2})-(\ref{d}),
we obtain 
\begin{eqnarray}\label{dH1}
-i\overline{H}_{\alpha }\omega _{12}\wedge \varphi +d\overline{H}_{\alpha
}\wedge \varphi &=&d\overline{H}_{\alpha -2}\wedge \omega _{\alpha -2,\alpha
}+d\overline{H}_{\alpha -1}\wedge \omega _{\alpha -1,\alpha }  \notag \\
&&+\omega _{\alpha -2,\alpha -1}\wedge \left( \overline{H}_{\alpha -2}\omega
_{\alpha -1,\alpha }-\overline{H}_{\alpha -1}\omega _{\alpha -2,\alpha
}\right) \\
&&+\overline{H}_{\alpha +1}\varphi \wedge \omega _{\alpha +1,\alpha }  \notag
\end{eqnarray}%
for $\alpha =2r+3,r\geq 0,$ and%
\begin{eqnarray}\label{dH2}
-i\overline{H}_{\alpha }\omega _{12}\wedge \varphi +d\overline{H}_{\alpha
}\wedge \varphi &=&d\overline{H}_{\alpha -3}\wedge \omega _{\alpha -3,\alpha
}+d\overline{H}_{\alpha -2}\wedge \omega _{\alpha -2,\alpha }  \notag \\
&&+\omega _{\alpha -3,\alpha -2}\wedge \left( \overline{H}_{\alpha -3}\omega
_{\alpha -2,\alpha }-\overline{H}_{\alpha -2}\omega _{\alpha -3,\alpha
}\right) \\
&&+\overline{H}_{\alpha -1}\varphi \wedge \omega _{\alpha -1,\alpha }  \notag
\end{eqnarray}%
for $\alpha =2r+4,r\geq 0.$

Each plane bundle $N^{r}f$ inherits a Riemannian connection from that of the
normal bundle  whose  curvature $K_{r}^{\ast}$  is given by 
$$d\omega_{2r+1,2r+2}=-K_{r}^{\ast }\omega_{1} \wedge \omega_{2}.$$

\begin{proposition}\cite{A}\label{Asp}
\textit{The  curvature $K_{r}^{\ast }$ of each plane bundle $N^{r}f$
of a minimal surface }$f:(M,ds^{2})\rightarrow \S^{n}$ \textit{is given by} 
\begin{equation*}
K_{1}^{\ast }=K_{1}^{\perp }-{\frac{\Vert B_{2}\Vert ^{2}}{2K_{1}^{\perp }}};%
{\ }K_{r}^{\ast }={\frac{}{}}{\frac{%
K_{r}^{\perp }\Vert B_{r-1}\Vert ^{2}}{2^{r-2}(K_{r-1}^{\perp })^{2}}}-{\frac{\Vert B_{r+1}\Vert ^{2}}{%
2^{r}K_{r}^{\perp }}},{\ }2\leq r\leq m-1.
\end{equation*}
\end{proposition}

We use the above mentioned notation throughout the paper.

\section{Hopf differential and the  $a$-invariants}
The Hopf differentials are
defined   from the higher fundamental forms and the complex structure in the following way. 
The complexified tangent bundle $TM\otimes \mathbb{C}$ is decomposed into
the eigenspaces of the complex structure $J$, called $T^{\prime }M$ and $%
T^{\prime \prime }M$, corresponding to the eigenvalues $i$ and $-i$. The
complex structure of $M$ is given by the orientation and the induced metric.
The $(r+1)$-th fundamental form $B_{r}$, which takes values in $N^{r}f$, can
be complex linearly extended to $TM\otimes \mathbb{C}$ with values in the
complexified vector bundle $N^{r}f\otimes \mathbb{C}$ and then decomposed
into its $(p,q)$-components, $p+q=r+1,$ which are tensor products of $p$ many 
1-forms vanishing on $T^{\prime \prime }M$ and $q$ many 1-forms vanishing on $%
T^{\prime }M$. The  minimality of $f$
implies that  the $(p,q)$-components of $B_{r}$
vanish, unless $p=r+1$ or $p=0.$ Hence, for a local complex
coordinate $z=x+iy$ on $M,$ we have the  decomposition 
\begin{equation*}
B_{r}=B_{r}^{(r+1,0)}dz^{r+1}+B_{r}^{(0,r+1)}d\bar{z}^{r+1},
\end{equation*}
where 
\begin{equation*}
B_{r}^{(r+1,0)}=B_{r}(\partial,...,\partial ){\ }{\ }{\hbox {and}}{\ }{\ }%
\partial ={\frac{1}{2}}\big({\frac{\partial }{\partial x}}-i{\frac{\partial 
}{\partial y}}\big ).
\end{equation*}%

The Hopf differentials are the differential forms%
\begin{equation}\label{H}
\Phi _{r}:=\langle B_{r}^{(r+1,0)},B_{r}^{(r+1,0)}\rangle dz^{2r+2},
\end{equation}%
of type $(2r+2,0),r=1,...,[(n-1)/2],$ where  $\langle.,.\rangle $ denotes the extension
of the usual Riemannian metric of $\S^{n}$\ to a complex-valued complex
bilinear form. These forms are defined at generic points and are independent
of the choice of coordinates, while $\Phi _{1}$ is globally well defined. It is a  consequence of the
structure equations that $\Phi _{1}$ is always holomorphic (cf. \cite{Ch}). 

Using (\ref{n1}) and (\ref{l}), we see that\textit{%
\begin{equation}\label{hl}
\left\vert \langle B_{r}^{(r+1,0)},B_{r}^{(r+1,0)}\rangle \right\vert ^{2}=%
\frac{F^{2r+2}}{2^{2r+4}}\left( \left\Vert B_{r}\right\Vert
^{4}-4^{r}(K_{r}^{\perp })^{2}\right).
\end{equation}%
}Hence, the zeros of $\Phi _{r}$ are precisely the points where $%
\mathcal{E}\sb r$ is a circle.

We choose a local complex coordinate $z=x+iy$ such that $\varphi =Fdz.$ 
Using the null basis $\eta_{r}=e_{2r+1}+ie_{2r+2},\overline{\eta}_{r}=e_{2r+1}-ie_{2r+2}$
 of the complexified bundle
$N^{r}f\otimes \mathbb{C}$ we have 
\begin{equation*}
 \langle B_{r}^{(r+1,0)},B_{r}^{(r+1,0)}\rangle =
\langle B_{r}^{(r+1,0)},\eta_{r}\rangle \langle B_{r}^{(r+1,0)},\overline{\eta}_{r}\rangle.
\end{equation*}
Then, from  (\ref{H}) we have
\begin{equation*}
4\Phi _{r}=\left( {\overline{H}_{2r+1}^{2}}+{\overline{H}%
_{2r+2}^{2}}\right) \varphi ^{2r+2}=k^{+}_{r}k^{-}_{r}\varphi^{2r+2},
\end{equation*}
where 
\begin{equation*}
k^{\pm}_{r}:= {\overline{H}_{2r+1}} \pm i{\overline{H}%
_{2r+2}}.
\end{equation*}

We introduce the $a$-\textit{invariants} as the functions
\begin{equation*}
a^{\pm}_{r}:= |k^{\pm}_{r}|.
\end{equation*}
They
 determine the geometry of the $r$-th curvature ellipse, since (\ref{n})-(\ref{l1}) imply that 
\begin{equation*}
a^{\pm}_{r}= \Big(2^{-r}\left\Vert B_{r}\right\Vert ^{2} \pm K_{r}^{\perp }\Big)^{1/2},
\end{equation*}
or equivalently
\begin{equation*}
a^{\pm}_{r}= \kappa _{r}{\pm}  \mu _{r},  
\end{equation*}
where $N^r f$ is equipped with orientation induced by the one of $M$ and $B_r$.  If
 $\Phi_{r}=0$, then  $a^{-}_{r}=0$ vanishes.

It will be convenient to set $a^{+}_{0}=2$ and $a^{-}_{0}=0$. 
We now recall the definition of exceptional surfaces.
\begin{definition}
A minimal surface $f:(M,ds^{2})\rightarrow \S^{n}$ is said to be exceptional if and only if all its Hopf differentials
 are holomorphic.
\end{definition}

The following provides a characterisation
 for exceptional surfaces  in terms of the higher curvature ellipses. 

\begin{theorem} \cite{V2}\label{exp}
A minimal surface is exceptional if and only if 
 its higher curvature ellipses
have constant eccentricity up to the last but one. 
\end{theorem}

The following  lemma is needed for the proof of the main results.

\begin{lemma}\label{dk}
For any exceptional surface  $f:(M,ds^{2})\rightarrow \S^{n}$ and  for any $1\leq r\leq m,$ where $m=[(n-1)/2]$, 
the following holds:
\begin{equation*}
 dk^{\pm}_{r}-i(s+1)k^{\pm}_{r}\omega_{12}{\pm}ik^{\pm}_{r}\omega_{2r+1,2r+2} \equiv 0{\ }{\rm {mod}}{\ }\varphi.
\end{equation*}
\end{lemma}
\begin{proof}
Using (\ref{dH1}), (\ref{dH2}) and arguing as in the proof of Proposition 4 in \cite{V2}, we obtain
\begin{equation*}
d\overline{H}_{2r+1}-i(r+1)\overline{H}_{2r+1}\omega _{12}-\overline{H}%
_{2r+2}\omega _{2r+1,2r+2}\equiv 0{\ }{\mbox {mod}}{\ }\overline{\varphi }
\end{equation*}%
and
\begin{equation*}
d\overline{H}_{2r+2}-i(r+1)\overline{H}_{2r+2}\omega _{12}+\overline{H}%
_{2r+1}\omega _{2r+1,2r+2}\equiv 0{\ }{\mbox {mod}}{\ }\overline{\varphi }.
\end{equation*}%
Then the lemma follows directly.
\end{proof}

\section{A characterisation of exceptional surfaces}
In this section, we give our main results according to which all exceptional surfaces are determined up 
to a multiparemeter family of isometric minimal deformations by the  $a$-invariants, 
provided that they satisfy certain restrictions.

For the proof of the results, we use the notion of absolute
value type functions introduced in \cite{EGT,ET}. A smooth complex valued function $p$ defined on a connected 
oriented surface $(M,ds^2)$  is 
called of \textit{%
holomorphic type}  if locally $p=p_{0}p_1,$ where $p_{0}$ is holomorphic and $p_{1}$ is  smooth without zeros.
  A function $a:M\rightarrow
\lbrack 0,+\infty )$  on  $M$ is called of \textit{%
absolute value type (AVT)}  if there is a function $p$  of holomorphic type on $M$  such that $a=|p|$. The zero set 
of such a function is either isolated
 or the whole of $M$, and outside its zeros the function is smooth. We need the following  lemmas that were proved in \cite{EGT,ET}.

\begin{lemma}\label{ht}
Let $p$ be a smooth complex valued function defined on $M, p\neq0$, and $\omega$ a real valued 1-form on $M$. Let $\psi:=pdz$
 for some conformal coordinate $z$. Then the equality
\begin{equation*}
 d\psi=i\omega\wedge\psi
\end{equation*} 
is valid if and only if $p$ is of holomorphic type and
\begin{equation*}
 \omega=2 {\rm{Im}}\big( \overline{\partial}(\log p)d\overline{z} \big).
\end{equation*} 
Moreover, 
 \begin{equation*}
d\omega=-\frac{1}{2i} \Delta \log |p| \overline{\varphi}\wedge\varphi.
\end{equation*} 
\end{lemma}

It is worth mentioning that $\overline{\partial}(\log p)$ and $\Delta \log |p|$ are well defined even at the zeros
 of $p$, if $p$ is of holomorphic type. Here, $\Delta$ denotes the Laplace-Beltrami operator of $(M,ds^2)$.

\begin{lemma}\label{log}
 Let $a$ be an AVT function on an open, simply connected subset $U$ of $\mathbb{C}$ such
 that $\varDelta^{0} \log a =0,$ where $\varDelta^{0}$ is the Euclidean Laplacian. Then there exists a
 holomorphic function $h$ on $U$ with $a=|h|.$
\end{lemma}

The function $h$ in Lemma \ref{log} is determined up to a factor $e^{i\theta}$.

\bigskip 
For any minimal surface we  set
$$\rho _{r}:=2^{r}\frac{|K_{r}^{\perp }|}{\left\Vert B_{r}\right\Vert ^{2}} {\ }{\ }{\mbox {and} }{\ }{\ }\rho _{0}:=1.$$
Obviously, $\rho _{r}=1$ precisely at points where $\Phi_r$ vanishes. We may now state the main results.

\begin{theorem}\label{main1}
Let $f:(M,ds^{2})\rightarrow \S^{n}$ be an exceptional surface with Gaussian curvature $K$ 
 and $m=[(n-1)/2].$ Then the functions $\rho _{r}$ are constant for any $1\leq r\leq m-1$, the $a$-invariants
 are AVT and  satisfy the following: 
\begin{equation}\label{as}
 a^{-}_{r}=\sigma_{r} a^{+}_{r}, 0\leq r\leq m-1, {\ } \sigma_{r}:=\sqrt{\frac{1-\rho _{r}}{1+\rho _{r}}}, a^{+}_{1}=
\sqrt{(1+\rho _{1})(1-K)},
\end{equation}
\begin{equation}\label{a}
\Delta \log a^{\pm}_{r}=(r+1)K\mp\Big( \frac{\rho _{r}b_{r}^{2} }{\rho^{2} _{r-1}b_{r-1}^{2}} -
 \frac{b_{r+1}^{2} }{\rho _{r}b_{r}^{2}}               \Big), {\ }0\leq r\leq m-1,
\end{equation}%
where 
$$b_{r}:=\sqrt{2^{r}/(1+\rho _{r})}a^{+}_{r}, {\ }0\leq r\leq m-1 {\ }{\ } {\mathrm{and}}{\ }{\ } b_m:=\|B_m\|.$$
Moreover
\begin{equation}\label{am}
\Delta \log a^{\pm}_{m} = (m+1)K \mp \frac{2^{m}K_{m}^{\perp}}{\rho^{2} _{m-1}b_{m-1}^{2}}. 
\end{equation}%

\end{theorem}

Theorem \ref{main1} shows that $a^{\pm}_{r}$ are intrinsic for any $ 1\leq r\leq m-1$ or equivalently $\|B_r\|$
 and $|K_r^{\perp}|$ are intrinsic. In the case where $f$ lies  in $\S^{2m+1}\subset \S^{2m+2}$, we have
$K^{\bot}_{m}=0$ and $ a^{+}_{m-1}=a^{-}_{m-1}$ and so Theorem \ref{main1} shows that all $a$-invariants are intrinsic.
 Furthermore, the metric $(a^{+}_{m})^{\frac{2}{m+1}}ds^{2}$ is flat. Obviously, if the $a$-invariants are constant, then 
 the surface is flat, 
 case that was studied by Kenmotsu  \cite{K1,K2}.

Any exceptional surface with non-vanishing first Hopf differential 
satisfies the Ricci condition, namely the metric $d\widehat{s}^{2}=\sqrt{1-K}ds^{2}$ is
flat away from points where $K=1$ or equivalently   $\Delta \log (1-K)=4K.$ This follows immediately from 
(4.1) and (\ref{a}) for $r=1.$

The converse of Theorem \ref{main1} is also true.

\begin{theorem}\label{main2}
Let $(M,ds^{2})$ be a simply connected two-dimensional Riemannian manifold with Gaussian curvature $K \lvertneqq 1$.
 Let $0<\rho _{r}\leq 1, 0\leq r\leq m-1,$ be constant numbers with $\rho _{0}= 1.$
Assume that there exist AVT functions $a^{\pm}_{r}, 1\leq r\leq m-1,$ 
 that satisfy (\ref{as}) and (\ref{a}), where 
$$b_{r}:=\sqrt{2^{r}/(1+\rho _{r})}a^{+}_{r}, {\ }0\leq r\leq m-1.$$ 
Let $K^{\perp}:M\rightarrow \mathbb{R}$ be a smooth function satisfying the inequality
 $|K^{\perp}|\leq 2^{-m}{b_{m}^{2}}.$ 
If the functions
\begin{equation*}
 a^{\pm}_{m}:=\sqrt{2^{-m}b_{m}^{2}\pm K^{\perp}}
\end{equation*} 
are AVT and satisfy (\ref{am}),
then for any $\theta \in \mathbb{R}^t$, 
 there exists a 
minimal surface $f_{\theta}:(M,ds^{2})\rightarrow \S^{n}$, with $m=[(n-1)/2]$
whose $a$-invariants are precisely the AVT functions $a^{\pm}_{r}, 1\leq r\leq m.$ Furthermore, $f$ is exceptional,  $t$ is the number of non-vanishing Hopf differentials,
 and any other minimal immersion of $(M,ds^{2})$ into $\S^n$ having the same $a$-invariants is congruent to some
 $f_{\theta}$.

\end{theorem}

Superconformal surfaces have vanishing Hopf differentials up to the last but one. This means 
that $a_{r}^{-}=0$ for any $0\leq r \leq m-1$, and so the following corollary follows immediately from Theorems 2 and 3.

\begin{corollary}\label{superconformal}
Let $f:(M,ds^{2})\rightarrow \S^{n}$ be a superconformal surface with Gaussian curvature $K$ 
 and set $m=[(n-1)/2].$  Then 
 the functions
 $a^{+}_{r}, 1\leq r\leq m-1,$ are AVT and satisfy the following: 
\begin{equation}\label{supr}
\Delta \log a^{+}_{r}=(r+1)K -  \frac{b_{r}^{2} }{b_{r-1}^{2}} +
 \frac{b_{r+1}^{2} }{b_{r}^{2}}, {\ }0\leq r\leq m-1,
\end{equation}%
where 
$$a^{+}_{1}=
\sqrt{2(1-K)}, {\ } b_{r}:=2^{\frac{r-1}{2}}a^{+}_{r}, {\ }0\leq r\leq m-1 {\ }{\ } {\mathrm{and}}{\ }{\ } b_m:=\|B_m\|.$$
Moreover
\begin{equation}\label{supm}
\Delta \log a^{\pm}_{m} = (m+1)K \mp \frac{2^{m}K_{m}^{\perp}}{b_{m-1}^{2}}. 
\end{equation}%
Conversely, let $(M,ds^{2})$ be a simply connected two-dimensional Riemannian manifold   with Gaussian curvature $K  \lvertneqq 1$.
 Assume  that there exist AVT functions $a^{+}_{r}, 1\leq r\leq m-1,$ 
and a non-negative function $b_{m}$  that fulfill (\ref{supr}). If for a given smooth function $K^{\perp}:M\rightarrow \mathbb{R}$ 
satisfying the inequality $|K^{\perp}|\leq 2^{-m}{b_{m}^{2}},$ 
the functions
\begin{equation*}
 a^{\pm}_{m}:=\sqrt{2^{-m}b_{m}^{2}\pm K^{\perp}}
\end{equation*} 
are AVT and satisfy (\ref{supm}),
then for any $\theta \in \mathbb{R}$ there exists a 
minimal surface $f_{\theta}:(M,ds^{2})\rightarrow \S^{n}$, with $m=[(n-1)/2]$,
whose $a$-invariants are precisely the AVT functions $a^{\pm}_{r}, 1\leq r\leq m,$ with
$a_{r}^{-}=0$ for any $0\leq r \leq m-1$.
 Furthermore, $f$ is superconformal
 and any other minimal immersion of $(M,ds^{2})$ into $\S^n$ having the same $a$-invariants is congruent to some $f_{\theta}$.
Furthermore, if $a_m^- = 0$, then
 $n$ is even and $f$ is rigid.
\end{corollary}

Miyaoka \cite{M1} determined all superconformal
surfaces lying  in  spheres of odd dimension in terms of 
solutions of the corresponding affine Toda equations. The above corollary gives another characterization of superconformal 
surfaces in any sphere.
 The class of superminimal surfaces has been
investigated by various authors (cf. \cite{B, Cal,Ch,J}). As a result,
superminimal surfaces are rigid, lie 
 in even dimensional
spheres. These results follow from Corollary \ref{superconformal}. 
Furthermore, Theorems 2 and 3 extend earlier results due to Eschenburg, Tribuzy and Guadalupe \cite{TG,ET0}.

\begin{proof}[Proof of Theorem \ref{main1}]
The fact that $\rho _{r}, 1\leq r\leq m-1,$ are constant follows immediately from Theorem \ref{exp}. Moreover, (4.1) is a consequence 
of the definition of $a$-invariants. 
 We choose the frame in the normal bundle as in Section 2 and we put
\begin{equation*}
 \psi_{r}^{\pm}:= k_{r}^{\pm}\varphi, {\ }     \omega_{r}^{\pm}:= r\omega_{12}\mp \omega_{2r+1,2r+2},   {\ }1\leq r\leq m.
\end{equation*}
Assume that $\varphi=\mu dz$ for a local conformal coordinate $z$ and set $\lambda =|\mu|.$ Appealing to Lemma \ref{dk}, we have
\begin{equation}
 d\psi_{r}^{\pm}=i\omega_{r}^{\pm}\wedge\psi_{r}^{\pm}, {\ }1\leq r\leq m.
\end{equation} 
Then  Lemma \ref{ht} implies that the functions
 $a^{\pm}_{r}, 1\leq r\leq m,$ are  AVT. Moreover,
 \begin{equation*}
d\omega_{r}^{\pm}=-\frac{1}{2i} \Delta \log |p_{r}^{\pm}| \overline{\varphi}\wedge\varphi,
\end{equation*} 
where $\psi_{r}^{\pm}=p_{r}^{\pm}dz=k_{r}^{\pm}\mu dz$, or  equivalently
 \begin{equation*}
rd\omega_{12}\mp d\omega_{2r+1,2r+2}=-\frac{1}{2i} \Delta \log (a_{r}^{\pm} \lambda )\overline{\varphi}\wedge\varphi.
\end{equation*} 
Now using (\ref{K}),  we obtain
\begin{equation*}
 \Delta \log (a_{r}^{\pm} \lambda)= rK \mp K^{*}_{r}.
\end{equation*}
Since $\Delta \log \lambda= -K $, we have 
\begin{equation*}
 \Delta \log a_{r}^{\pm}= (r+1)K \mp K^{*}_{r}, {\ } 1\leq r \leq m.
\end{equation*}
Thus (\ref{a}) and (\ref{am}) follow from this, Proposition \ref{Asp} and the equations 
$$(a_{r}^{+})^{2}=
2^{-r}(1+\rho _{r})\left\Vert B_{r}\right\Vert ^{2}=\dfrac{1+\rho _{r}}{\rho _{r}}K_{r}^{\perp}.$$
\end{proof}

\begin{proof}[Proof of Theorem \ref{main2}]
 \textit{(i) Existence}. Let $U\subset M$ be an open, simply connected subset and $z=x+iy$ a conformal coordinate on $U$.
Choose an orthonormal frame $\{e_{1},e_{2}\}$ so that $\varphi =\mu dz$ and set $\lambda =|\mu|$.

Since  the functions  $a^{\pm}_{r}, 1\leq r\leq m-1,$ are AVT, there exist functions 
$k^{\pm}_{r}:U\rightarrow \mathbb{C}$ of holomorphic type such that $a^{\pm}_{r}=|k^{\pm}_{r}|$. In particular, we set $k^{+}_{0}=2$
and $k^{-}_{0}=0.$ We put 
\begin{equation*}
 \psi_{r}^{\pm}:= k_{r}^{\pm}\varphi=p_{r}^{\pm}dz,  
  {\ }1\leq r\leq m.
\end{equation*}
Appealing to Lemma \ref{ht}, we have 
\begin{equation}\label{ps}
 d\psi_{r}^{\pm}=i\omega_{r}^{\pm}\wedge\psi_{r}^{\pm}, {\ } 1\leq r\leq m,
\end{equation} 
where 
\begin{equation}\label{omr}
 \omega_{r}^{\pm}=2 {\rm{Im}}\big( \overline{\partial}(\log p_{r}^{\pm})d\overline{z} \big).
\end{equation} 
Moreover, 
 \begin{equation}\label{domr}
d\omega_{r}^{\pm}=-\frac{1}{2i} \Delta \log |p_{r}^{\pm}| \overline{\varphi}\wedge\varphi.
\end{equation} 

For any  $\theta
_{r} \in \mathbb{R}, 0\leq r\leq m-2$, we define the  forms
\begin{eqnarray*}
{\omega }_{2r+1,2r+3} &:=&\frac{(1+\rho _{r})(1+\sigma_{r+1})(1-\sigma _{r})}{2\rho _{r}}  \mathrm{Re}\big(\frac{e^{i\theta
_{r+1}}k^{+}_{r+1}}{e^{i\theta
_{r}}k^{+}_{r}} \varphi\big),\\
{\omega }_{2r+2,2r+3} &:=&-\frac{(1+\rho _{r})(1+\sigma_{r+1})(1+\sigma _{r})}{2\rho _{r}}  \mathrm{Im}\big(\frac{e^{i\theta
_{r+1}}k^{+}_{r+1}}{e^{i\theta
_{r}}k^{+}_{r}} \varphi\big), \\
{\omega }_{2r+1,2r+4} &:=&\frac{(1+\rho _{r})(1-\sigma_{r+1})(1-\sigma _{r})}{2\rho _{r}}  \mathrm{Im}\big(\frac{e^{i\theta
_{r+1}}k^{+}_{r+1}}{e^{i\theta
_{r}}k^{+}_{r}} \varphi\big), \\
{\omega }_{2r+2,2r+4} &:=&\frac{(1+\rho _{r})(1-\sigma_{r+1})(1+\sigma _{r})}{2\rho _{r}}  \mathrm{Re}\big(\frac{e^{i\theta
_{r+1}}k^{+}_{r+1}}{e^{i\theta
_{r}}k^{+}_{r}} \varphi\big),
\end{eqnarray*}%
whereas 
for any $1\leq r\leq m-1,$ we define
\begin{equation*}
{\omega }_{2r+1,2r+2}:=\left\{ 
\begin{array}{c}
r{\omega }_{12}-\omega_{r}^{+}{\ }{\ } \mathrm{if} {\ } \rho _{r}=1\\
 
\frac{1}{2}(\omega_{r}^{+}-\omega_{r}^{-}) {\ }{\ } \mathrm{if} {\ } \rho _{r}<1.\\ %
\end{array}%
\right.
\end{equation*}%
Furthermore,  for any $\theta
_{m} \in \mathbb{R}$,  we define
\begin{eqnarray*}
{\omega }_{2m-1,2m+1} &:=&\frac{(1+\rho _{m-1})(1-\sigma_{m-1})}{2\rho _{m-1}}  \mathrm{Re}\big(\frac{e^{i\theta
_{m}}(k^{+}_{m}+k^{-}_{m})}{e^{i\theta
_{m-1}}k^{+}_{m-1}} \varphi\big),\\
{\omega }_{2m,2m+1} &:=&-\frac{(1+\rho _{m-1})(1+\sigma_{m-1})}{2\rho _{m-1}}  \mathrm{Im}\big(\frac{e^{i\theta
_{m}}(k^{+}_{m}+k^{-}_{m})}{e^{i\theta
_{m-1}}k^{+}_{m-1}} \varphi\big), \\
{\omega }_{2m-1,2m+2} &:=&\frac{(1+\rho _{m-1})(1-\sigma_{m-1})}{2\rho _{m-1}}  \mathrm{Im}\big(\frac{e^{i\theta
_{m}}(k^{+}_{m}-k^{-}_{m})}{e^{i\theta
_{m-1}}k^{+}_{m-1}} \varphi\big), \\
{\omega }_{2m,2m+2} &:=&\frac{(1+\rho _{m-1})(1+\sigma _{m-1})}{2\rho _{m-1}}  \mathrm{Re}\big(\frac{e^{i\theta
_{m}}(k^{+}_{m}-k^{-}_{m})}{e^{i\theta
_{m-1}}k^{+}_{m-1}} \varphi\big),\\
{\omega }_{2m+1,2m+2}&:=&\frac{1}{2}(\omega_{m}^{+}-\omega_{m}^{-}).
\end{eqnarray*}%
 In all other cases, we set ${\omega }_{AB}=0.$ 

Our aim is to prove that the forms $
{\omega }_{j}$ and ${\omega }_{AB}$ satisfy the structure equations (\ref{s1}) and (\ref{s2}). We will only confirm that
\begin{equation}\label{dom1}
d\omega _{2r+1,2r+2} =\sum_{C}\omega _{2r+1,C}\wedge \omega _{C,2r+2}, {\ } 1\leq r\leq m
\end{equation} 
and 
\begin{equation}\label{dom2}
d\omega _{2r+1,2r+3} =\sum_{C}\omega _{2r+1,C}\wedge \omega _{C,2r+3}, {\ } 1\leq r\leq m.
\end{equation} 
The proof of the rest structure equations follows in the same manner. 

From (\ref{domr}), our assumption (\ref{a}) and $\Delta \log \lambda= -K $, we have
\begin{equation*}
2id \omega ^{\pm}_{r}=\Big\{-rK\pm      \Big( \frac{\rho _{r}b_{r}^{2} }{\rho^{2} _{r-1}b_{r-1}^{2}} -
 \frac{b_{r+1}^{2} }{\rho _{r}b_{r}^{2}}  \Big)  \Big\}\overline{\varphi}\wedge\varphi, {\ } 1\leq r\leq m-1.
\end{equation*}
If $\rho_{r}=1,$ then ${\omega }_{2r+1,2r+2}=r{\omega }_{12}-\omega_{r}^{+} $ and 
on account of (\ref{K}), we find
\begin{equation}\label{dom3}
  d\omega _{2r+1,2r+2} =-\dfrac{1}{2i} \Big( \frac{\rho _{r}b_{r}^{2} }{\rho^{2} _{r-1}b_{r-1}^{2}} -
 \frac{b_{r+1}^{2} }{\rho _{r}b_{r}^{2}}  \Big)   \overline{\varphi}\wedge\varphi, {\ } 1\leq r\leq m-1.
\end{equation}
If $\rho_{r}<1,$ then ${\omega }_{2r+1,2r+2}=\frac{1}{2}(\omega_{r}^{+}-\omega_{r}^{-})$. From (\ref{domr}) and (\ref{as}), we obtain
\begin{equation*}
  d\omega _{2r+1,2r+2} =0, {\ } 1\leq r\leq m-1.
\end{equation*}
Moreover, from (\ref{a}) we have 
\begin{equation*}
  \frac{\rho _{r}b_{r}^{2} }{\rho^{2} _{r-1}b_{r-1}^{2}} -
 \frac{b_{r+1}^{2} }{\rho _{r}b_{r}^{2}} =0.
\end{equation*}
This shows that  (\ref{dom3}) is true in either case. On the other hand, by the definition of the forms $\omega_{AB}$, 
we see that
\begin{equation*}
 \sum_{C}\omega _{2r+1,C}\wedge \omega _{C,2r+2}=-\dfrac{i}{2}\Big( \frac{b_{r+1}^{2}}{\rho_{r}b_{r}^{2}} -
 \frac{\rho_{r}b_{r}^{2} }{\rho _{r-1}^{2}b_{r-1}^{2}}  \Big)   \overline{\varphi}\wedge\varphi, {\ } 1\leq r\leq m-1.
\end{equation*}
Hence, (\ref{dom1}) is satisfied for any $1\leq r\leq m-1$. 
Using (\ref{domr}), the assumption (\ref{a}), and ${\omega }_{2m+1,2m+2}=\frac{1}{2}(\omega_{m}^{+}-\omega_{m}^{-})$, we find 
\begin{equation*}
 d\omega _{2m+1,2m+2}=-\dfrac{1}{2i}\dfrac{2^{m}K^{\bot}}{\rho^{2}_{m-1}b^{2}_{m-1}}\overline{\varphi}\wedge\varphi.
\end{equation*}
By direct calculations we conclude  that
\begin{equation*}
 \sum_{C}\omega _{2m+1,C}\wedge \omega _{C,2m+2}=
\dfrac{i}{2}\dfrac{2^{m}K^{\bot}}{\rho^{2}_{m-1}b^{2}_{m-1}}\overline{\varphi}\wedge\varphi.
\end{equation*}
Consequently, (\ref{dom1}) holds true for any $1\leq r\leq m$.

Now from (\ref{ps}), $\psi ^{\pm}_{r}=k^{\pm}_{r}\varphi$ and (\ref{d}), we have 
\begin{equation}\label{dkr1}
 dk^{\pm}_{r}\wedge\varphi-i(r+1)k^{\pm}_{r}(\omega_{12}+{\omega}^{\pm}_{r})\wedge
\varphi=0, {\ } 1\leq r\leq m.
\end{equation}
We claim that
\begin{equation}\label{dkr2}
 dk^{\pm}_{r}\wedge\varphi-i(r+1)k^{\pm}_{r}\omega_{12}\wedge\varphi \pm ik^{\pm}_{r}\omega_{2r+1,2r+2}\wedge\varphi=0,
 {\ } 1\leq r\leq m.
\end{equation} 
Obviously, (\ref{dkr1}) easily implies (\ref{dkr2})  if $\rho_{r}=1$ and ${\ } 1\leq r\leq m-1$,
since $\omega^{+}_{r}=r\omega_{12}-\omega_{2r+1,2r+2}$. 

Assume that  $\rho_{r}<1$ and ${\ } 1\leq r\leq m-1$. Then our assumptions (\ref{as}) and (\ref{a}) yield $\Delta \log a_{r}^{+}=(r+1)K$ and 
since $\Delta \log \lambda =-K$, we have $\Delta \log (a_{r}^{+}\lambda^{r+1})=0$. 
 According to Lemma \ref{log}, we deduce that there exists a holomorphic function $g_{r}^{+}$ such
 that $a_{r}^{+}\lambda^{r+1}=|g_{r}^{+}|$. Moreover, $a_{r}^{-}\lambda^{r+1}=|g_{r}^{-}|$, where $g_{r}^{-}:=\sigma_{r}g_{r}^{+}$.

We may choose $k^{\pm}_{r}$ so that $k^{+}_{r}=\mu^{-r-1}g_{r}^{+}$ and $k^{-}_{r}=\sigma_{r}k_{r}^{
+}$. Since $k^{+}_{r}k_{r}^{-}\mu^{2r
+2}$ is holomorphic, from (\ref{omr}) we find
$$-\dfrac{1}{2r}(\omega^{+}_{r}+\omega^{-}_{r})=2 {\rm{Im}}\big( \overline{\partial}(\log \mu)d\overline{z} \big),$$
or equivalently, in view of $\varphi=\mu dz$,
$$d\varphi=-i\dfrac{\omega^{+}_{r}+\omega^{-}_{r}}{2^{r}}\wedge\varphi.$$
From (\ref{d}), we deduce that
$$\omega_{12}=\dfrac{1}{2r}(\omega^{+}_{r}+\omega^{-}_{r})$$
and so 
$$\omega^{\pm}_{r}=r \omega_{12}\mp \omega_{2r+1,2r+2}.$$
Thus from (\ref{ps}) and (2.5), we infer that (\ref{dkr2}) holds for any $1\leq r\leq m-1$. We note that since $k^{-}_{r}=\sigma_{r}k_{r}^{
+}$, (\ref{dkr2}) yields $ \omega_{2r+1,2r+2}=0$ for any $1\leq r\leq m-1$, with $\rho_{r}<1$.

It remains to prove  (\ref{dkr2}) for $r=m$. 
 Our assumption (\ref{am}) yields 
$$\Delta \log (a_{m}^{+}a_{m}^{-})=2(m+1)K.$$
 Using
 $\Delta \log \lambda =-K$, we find $\Delta \log (a_{m}^{+}a_{m}^{-}\lambda^{2m+2})=0$. 
According to Lemma \ref{log}, there exists a holomorphic function $g_{m}$ such
 that 
$$a_{m}^{+}a_{m}^{-}\lambda^{2m+2}=|g_{m}|.$$
 We may choose $k^{\pm}_{m}$ so that $k^{+}_{m}k^{-}_{m}\mu^{2m+2}=g_{m}$. 
Since  $k_{m}^{+}k^{-}_{m}\mu^{2m+2}$ is holomorphic, from (\ref{omr}) we have 
$$-\dfrac{1}{2m}(\omega^{+}_{m}+\omega^{-}_{m})=2 {\rm{Im}}\big( \overline{\partial}(\log \mu)d\overline{z} \big).$$
We may now argue as above to deduce that (\ref{dkr2}) holds for $r=m$.

The exterior derivative of $\omega_{2r+1,2r+2}$, for any $1\leq r\leq m$, is computed by using (\ref{dkr2}) and (2.5). On the other hand, since $\omega_{2r+1,2r+2}=0$ for any $1\leq r\leq m-1$ with $\rho_{r}<1$, from the definition of the forms ${\omega }_{AB},$
   we easily see that (\ref{dom2}) holds true.

Now the fundamental theorem of
submanifolds implies that there exists an isometric immersion $f_{\theta
_{1},...,\theta _{m}}:(U,ds^{2})\rightarrow \S^{n}$ with
corresponding connection forms ${\omega }_{AB}.$ 
Clearly,  ${f}$ is minimal and the components of its complexified higher fundamental forms are given by 
\begin{equation*}
{H}_{2r+1}{\omega }_{2r+1,\alpha }+{H}_{2r+2}%
{\omega }_{2r+2,\alpha }={H}_{\alpha }\overline{\varphi }%
\text{ \ for }\alpha =2r+3,2r+4,
\end{equation*}%
$0\leq r\leq m-1$, where ${H}_{1}=1,{H}_{2}=i,$
and when $n=2m+1$ 
\begin{equation*}
{H}_{2m-1}{\omega }_{2m-1,2m+1}+{H}_{2m}%
{\omega }_{2m,2m+1}={H}_{2m+1}\overline{\varphi }.
\end{equation*}%
By induction and by the definition of the forms $\omega_{AB}$, we deduce that 
\begin{equation*}
 \overline{H}_{2r+1}=\dfrac{1}{2}e^{i\theta_{r}}(1+\sigma_{r})k_{r}^{+} {\ }{\ } \mathrm{and} {\ }{\ }\overline{H}_{2r+2}
=-\dfrac{i}{2}e^{i\theta_{r}}(1-\sigma_{r})k_{r}^{-}, {\ }1\leq r\leq m-1, 
\end{equation*} 
whereas
\begin{equation*}
 \overline{H}_{2m+1}=\dfrac{1}{2}e^{i\theta_{m}}(k_{m}^{+}+k_{m}^{-}) {\ }{\ } \mathrm{and} {\ }{\ }\overline{H}_{2m+2}
=-\dfrac{i}{2}e^{i\theta_{m}}(k_{m}^{+}-k_{m}^{-}). 
\end{equation*} 
Since $k_{r}^{-}=\sigma_{r}k_{r}^{+}$ for any $1\leq r\leq m-1$,  we see that the Hopf 
differentials of $f$ are given by
$$\Phi_{r}=\dfrac{1}{4}(\overline{H}_{2r+1}^{2}+\overline{H}_{2r+2}^{2})\varphi^{2r+2}
=\dfrac{1}{4}e^{i\theta_{r}}k_{r}^{+}k_{r}^{-}\mu^{2r+2}dz^{2r+2}.$$
In particular, they are holomorphic and $f$ is exceptional. If $\rho_{r}=1$ for some $r$, then $\Phi_{r}=0$ and
 by  \cite{V}, $f$ does not depend on $\theta_{r}$.

To prove that $f$ is well defined on the whole of $M$, we cover $M$ with simply connected coordinate neighbourhoods $U_{a}$. 
Then we have exceptional surfaces $f_a:U_a\rightarrow \S^n$ which can be chosen so that they have the same
 Hopf differentials in the intersections $U_{a}\cap U_{b}$. This is achievable, since $M$ is simply connected. Thus by \cite{V}, 
$f_a$ and  $f_b$ are congruent on $U_{a}\cap U_{b}$. Continuing in this way, we obtain a minimal
 surface  $f:M\rightarrow \S^n$ that
has the desired properties.

\textit{(ii) Uniqueness.} Now let $\tilde{f}:(M, ds^2)\rightarrow \S^n$ be another minimal surface arising from the same data 
$a_{r}^{\pm}, 1\leq r\leq m$.
Then $f=f_{\theta
_{1},...,\theta _{m}}$ and $\tilde{f}$ have congruent higher curvature ellipses.
 By Theorem \ref{exp}, $\tilde{f}$ is also exceptional. 
Moreover, the Hopf differentials
 have the same length and are holomorphic. Hence there exist real numbers $\eta_{1},...,\eta_{m}$ so
 that $\widetilde{\Phi}_{r}=e^{i\eta _r}{\Phi}_{r}$ for any $r$. This means that $\tilde{f}$ and $f_{\theta
_{1}+\eta_{1},...,\theta _{m}+\eta_{m}}$ are congruent, since they have the same Hopf differentials (cf. \cite{V}).
\end{proof}

\section{Polar and self-dual surfaces}
The Gauss map of a minimal surface $M$ in $\S^3$ defines another minimal surface,
 the polar \cite{L0}   of $M$, which is conformal to $M$.  More generally, the \textit{polar} \cite{DF}
of a minimal surface 
 $f:(M,ds^{2})\rightarrow \S^{2m+1}$  is the map $f^{*}:M^* \rightarrow \S^{2m+1}$ defined by $f^{*}=e_{2m+1},$  where
$M^*$ is the set of generic points, and $e_{2m+1}$ is a unit section of the last normal bundle. If $f$ is exceptional, then 
$f^*$ is defined all over $M$, since all higher normal bundles are well 
defined over singular points \cite[Prop. 4]{V2}.
The polar of a superconformal surface is again a superconformal surface (cf. \cite{M1}).

\begin{proposition}\label{polar}
Let $f:(M,ds^{2})\rightarrow \S^{2m+1}$ be a  superconformal surface. Then $f$ and $f^*$ have the same Hopf differentials and the induced metric of $f^*$ is 
 $ds_*^2=\big({2a_{m}^{+}}/{a^{+}_{m-1}}\big)^{2}ds^{2}$. Moreover, 
 the polar of $f^*$ is $f$.
\end{proposition}
\begin{proof}
 According to 
 \cite[Prop. 8]{DF}, 
the polar is an elliptic surface (see \cite{DF}) whose higher normal bundles are given by $N^{r}{f^*}=N^{m-1-r}f$ 
for any $0\leq r\leq m-1$, where $N^{0}f:=df(TM^2)$ 
and $N^{0}f^{*}:=df^{*}(TM^2)$. Furthermore, $N^{m}{f^*}=\mathrm{span}\{ f^*\}$. Since $f$ is superconformal, all its
higher curvature ellipses are circles up to the last but one. Thus the corresponding complex structures $J_r$ and 
$\widetilde{J}_{r}=J_{m-1-r}^{t}$ defined in \cite{DF} are orthogonal. This means that $f^*$ is minimal, all its
higher curvature ellipses are circles up to the last but one, and so it is superconformal.  
 
To compute the last Hopf differential $\Phi_{m}^{*}$, we choose a local conformal coordinate $z$ and proceed as follows:
\begin{eqnarray*}
\Phi_{m}^{*}\!\!\!&=&\!\!\!\langle B_{m}^{*(m+1,0)},B_{m}^{*(m+1,0)}\rangle dz^{2m+2}
=\langle \overline{\nabla }_{\partial} \dots \overline{\nabla }_{\partial}df^{*}(\partial), f\rangle^{2}dz^{2m+2}\\
\!\!\!&=&\!\!\!\Big( \partial \langle \overline{\nabla }_{\partial} \dots \overline{\nabla }_{\partial}df^{*}(\partial), f\rangle-
 \langle \overline{\nabla }_{\partial} \dots \overline{\nabla }_{\partial}df^{*}(\partial), df(\partial)\rangle\Big)^{2}dz^{2m+2}\\
\!\!\!&=&\!\!\!\langle \overline{\nabla }_{\partial}  \dots \overline{\nabla }_{\partial}df^{*}(\partial), df(\partial)\rangle ^{2}dz^{2m+2}\\
\!\!\!&=&\!\!\!
\langle df^{*}(\partial), \overline{\nabla }_{\partial} \dots \overline{\nabla }_{\partial}df(\partial)\rangle ^{2}dz^{2m+2}\\
\!\!\!&=&\!\!\!\langle f^{*}, \overline{\nabla }_{\partial} \dots\overline{\nabla }_{\partial}df(\partial)\rangle ^{2}dz^{2m+2}\\
\!\!\!&=&\!\!\!\langle e_{2m+1}, \overline{\nabla }_{\partial} \dots\overline{\nabla }_{\partial}df(\partial)\rangle ^{2}dz^{2m+2}\\
\!\!\!&=&\!\!\!\langle B_{m}^{(m+1,0)},B_{m}^{(m+1,0)}\rangle dz^{2m+2}=\Phi_{m}.
\end{eqnarray*}
Since $f$ is superconformal, we may choose the frame so that  $H_{2r+1}=\kappa _{r}$ and $H_{2r+2}=i\kappa _{r}$ for
 any $1\leq r\leq m-1$. Then (\ref{01}) and (\ref{02}) yield
$$ \omega_{2m-1,2m+1}=\dfrac{1}{2\kappa _{m-1}}\big(H_{2m+1}\overline{\varphi}+  \overline{H}_{2m+1}\varphi  \big),$$
$$ \omega_{2m,2m+1}=\dfrac{i}{2\kappa _{m-1}}\big(H_{2m+1}\varphi- {H}_{2m+1}\overline{\varphi } \big).$$
Using (\ref{om}) and the above, for any $X,Y$ tangent to $M^2$, we find
$$\langle df^{*}(X), df^{*}(Y)\rangle=\dfrac{\kappa^2 _{m}}{\kappa^2 _{m-1}}\langle X,Y\rangle
=\Big(\dfrac{2a_{m}^{+}}{a^{+}_{m-1}}\Big)^{2}\langle X,Y\rangle.$$
\end{proof}
 
The polar of a superconformal surface $f:(M,ds^{2})\rightarrow \S^{2m+1}$ 
 can be characterised, up to congruence, as the minimal surface with
 the same Hopf differentials as $f$ and induced metric $\big(2{a_{m}^{+}}/{a^{+}_{m-1}}\big)^{2}ds^{2}$.

 The  following  is a consequence of Proposition \ref{polar}
 and the main result in \cite{V}.

\begin{proposition} \cite{M1}\label{isom}
 Let $f:(M,ds^{2})\rightarrow \S^{2m+1}$ be a  superconformal surface. Then $f$ and $f^*$ are isometric 
 if and only if they are congruent.
\end{proposition}

We are interested in superconformal surfaces  with the property that they are isometric to their polar.
 We call these surfaces \textit{self-dual}.

The following corollary provides a link between the Ricci condition and furnishes  examples of 
 self-dual surfaces.

\begin{corollary}
(i) A non-flat superconformal  surface in $\S^{2m+1}$ that satisfies the Ricci condition
is self-dual if and only if $m=4k+3$. 

(ii) Superconformal surfaces  in $\S^{4k+1}$ are self-dual if  they are flat.
\end{corollary}
\begin{proof}
 Part (i) follows from Lemma \ref{v1} and Proposition \ref{polar}, while part (ii) follows immediately
from Corollary \ref{superconformal} and Proposition \ref{polar}.
\end{proof}

We now give a  characterisation of self-dual surfaces. It turns out that this property is intrinsic. 
\begin{proposition}
 Let $f:(M,ds^{2})\rightarrow \S^{2m+1}$ be a  superconformal surface. If $f$ is self-dual, then its
$a$-invariants  satisfy
\begin{equation}\label{ratioa}
\dfrac{a^{+}_{m-r}}{{a^{+}_{m-r-1}}}=\dfrac{a^{+}_{r}}{{a^{+}_{r-1}}}, {\ } 0\leq r\leq m,
\end{equation}
where $a^{+}_{0}=:2, a^{+}_{-1}:=4$.  Moreover:

(i) If $m=2l$, then $a^{+}_{m}=\frac{1}{4}a^{+}_{l-1}a^{+}_{l}$ and $\Delta \log(a^{+}_{l-1}a^{+}_{l})=(m+1)K.$ 

(ii) If $m=2l+1$, then $a^{+}_{m}={\frac{1}{4}}(a^{+}_{l})^{2}$ and
$\Delta \log a^{+}_{l}=(m+1)K/2$. 

 Conversely, if $a^{+}_{m}=\frac{1}{4}a^{+}_{l-1}a^{+}_{l}$   when $m=2l$, or  $a^{+}_{m}={\frac{1}{4}}(a^{+}_{l})^{2}$ when 
 $m=2l+1$, then $f$ is self-dual.
\end{proposition}
\begin{proof}
 Assume that $f$ is self-dual. Proposition \ref{polar} shows that (\ref{ratioa}) holds for $r=0$. From this we obtain 
$\Delta \log a^+_m = \Delta \log a^+_{m-1},$ and appealing to Corollary \ref{superconformal}, we see that (\ref{ratioa}) holds for $r=1$. By reduction and  Corollary \ref{superconformal}, we prove (\ref{ratioa}) for any $ 0\leq r\leq m$. If $m=2l$, then (\ref{ratioa}) yields
 $a^{+}_{m}=\frac{1}{4}a^{+}_{l-1}a^{+}_{l}$ and consequently $\Delta \log(a^{+}_{l-1}a^{+}_{l})=(m+1)K$ follows from Corollary \ref{superconformal}.
The case  $m=2l+1$ is similar.

Conversely, we assume that $a^{+}_{m}=\frac{1}{4}a^{+}_{l-1}a^{+}_{l}$   and $m=2l$. The  case where $m=2l+1$ is treated in a similar manner. 
Using (\ref{supr}) and (\ref{supm}), we have 
$$\dfrac{a^{+}_{l+1}}{{a^{+}_{l}}}=\dfrac{a^{+}_{l-1}}{{a^{+}_{l-2}}}.$$
From this and  Corollary \ref{superconformal}, we obtain
$$\dfrac{a^{+}_{l+2}}{{a^{+}_{l+1}}}=\dfrac{a^{+}_{l}}{{a^{+}_{l-1}}}.$$
Inductively, we have that (\ref{ratioa}) holds for any $ 0\leq r\leq m$. In particular, this yields  $2{a_{m}^{+}}={a^{+}_{m-1}},$
and so by Proposition \ref{polar}, $f$ is self-dual.
\end{proof}

Now we characterize all metrics which arise as induced metrics on self-dual surfaces.

\begin{theorem}\label{sd}
Let $(M,ds^{2})$ be a simply connected two-dimensional Riemannian manifold   with Gaussian curvature $K \lvertneqq 1$.
 We consider the non-negative functions $a^{+}_{r}, 0\leq r\leq l,$ defined inductively by
 \begin{equation}\label{sda}
\Delta \log a^{+}_{r}=(r+1)K -  2\frac{a_{r}^{2} }{a_{r-1}^{2}} +2
 \frac{a_{r+1}^{2} }{a_{r}^{2}}, {\ }0\leq r\leq l-1,
\end{equation}%
where $a^{+}_{0}=:2, a^{+}_{-1}:=4, a^{+}_{1}:=\sqrt{2(1-K)}$ and $l$ being a positive integer.
Assume that these functions are AVT and  either $\Delta \log (a^{+}_{l-1}a^{+}_{l})=(2l+1)K$ or
 $\Delta \log a^{+}_{l}=(l+1)K$. Then for any $\theta \in \mathbb{R}$ there exists a self-dual 
surface $f_{\theta}:(M,ds^{2})\rightarrow \S^{n}$, with $n=4l+1$ or  $n=4l+3$, 
whose $a$-invariants up to order $l$ are precisely the AVT functions $a^{+}_{r}$, $1\leq r \leq l$.
 Furthermore, any other self-dual immersion of $(M,ds^{2})$ arising from these data  is congruent to some $f_{\theta}$. 
\end{theorem}
\begin{proof}
 Assume that $\Delta \log (a^{+}_{l-1}a^{+}_{l})=(2l+1)K$. We define  the AVT functions
 $a_{r}^{+}, l+1\leq r\leq 2l$, by
\begin{equation*}
\dfrac{a^{+}_{2l-r}}{{a^{+}_{2l-r-1}}}=\dfrac{a^{+}_{r}}{{a^{+}_{r-1}}}.
\end{equation*}
 Using (\ref{sda}) and 
$$\Delta \log (a^{+}_{l-1}a^{+}_{l})=(2l+1)K,$$
we 
can prove by induction that  $a_{r}^{+}, 1\leq r\leq 2l$, satisfy 
(\ref{supr}) and (\ref{supm}).  Corollary \ref{superconformal} implies that for any $\theta \in \mathbb{R}$ there exists a superconformal 
surface $f_{\theta}:(M,ds^{2})\rightarrow \S^{n}$, with $n=4l+1$. Obviously $f$ is self-dual.
 Furthermore, any other self-dual immersion of $(M,ds^{2})$ arising from these data  is congruent to some $f_{\theta}$.

The case where  $\Delta \log a^{+}_{l}=(l+1)K$ is treated in a similar manner, ending up with  superconformal self-dual
surfaces in $ \S^{4l+3}$. 
\end{proof}

The following provides a characterisation of self-dual surfaces in $\S^5$.

\begin{corollary}
 For any  superconformal surface $f:(M,ds^{2})\rightarrow \S^{5}$  the following are equivalent:

(i) $f$ is self-dual.

(ii) $f$  and $f^*$ have the same Gaussian curvature.

(iii) $f$ is  congruent to pseudoholomorphic curve that lies in $\S^5$.
\end{corollary}
\begin{proof}
Proposition \ref{polar} implies  that the Gaussian curvature of  $f^{*}$
is  
$$K_{*}=\Big(\frac{a^{+}_{1}}{2a_{2}^{+}}\Big)^{2} \Big(K-\Delta \log {a^{+}_{2}}+  \Delta \log {a^{+}_{1}}  \Big).$$
Appealing to Corollary \ref{superconformal}, we find that 
$$K_{*}=1-\dfrac{(a^{+}_{1})^4}{8(a_{2}^{+})^2}.$$
From $K=1-{(a^{+}_{1})^2}/{2}$
 and the above, we deduce that  $K_{*}=K$ is equivalent to  $2a_{2}^{+}=a^{+}_{1}$, which by virtue of Proposition \ref{polar} shows 
 the equivalence between  (i) and (ii). The equivalence between 
(i) and (iii) follows from Theorem 4(ii).
\end{proof}

\section{Pseudoholomorphic curves in $\S^{6}$}
It is well known that the multiplicative structure on the Cayley numbers $\mathbb{O}$ can be used to define an
 almost complex structure $J$ on the sphere $\S^{6}$ in $\mathbb{R}^7$. This complex structure is not integrable but is
 nearly K\"{a}hler. A pseudoholomorphic curve \cite{B}  is a non-constant map $f:M\rightarrow \S^{6}$ whose
 differential  is complex linear,  $M$ being a Riemann surface. 
Pseudoholomorphic curves in $\S^{6}$ are  superconformal.

As an application of  the main  results ,  we provide both
an extrinsic and an intrinsic characterisation for each type  (\cite{BVW}) of pseudoholomorphic curves, via the $a$-invariants,  among the class of superconformal  
 surfaces. An intrinsic characterisation was given by Hashimoto \cite{H}. For another approach see \cite{EV}. 

\begin{theorem}\label{pseudo}
 Let $f:(M,ds^{2})\rightarrow \S^{6}$ be a superconformal surface with Gaussian curvature $K$. Then the following hold:

(i) $f$ is locally $O(7)$-congruent to a superminimal pseudoholomorphic curve if and only if $\Delta \log (1-K)=6K-1,$
 or equivalently if and only if $f$ is superminimal with $a_{2}^{+}=a_{1}^{+}/2$.

(ii) $f$ is locally $O(7)$-congruent to a pseudoholomorphic curve which lies in $\S^5$ if and only if $\Delta \log (1-K)=6K,$
 or equivalently if and only if $a_{2}^{+}=a_{2}^{-}=a_{1}^{+}/2$ and $a_{1}^{-}=0$.

(iii) $f$ is locally $O(7)$-congruent to a non-superminimal  pseudoholomorphic curve if and only if $6K>\Delta \log (1-K)>6K-1$ and 
$$\Delta \log \Big(\big(1-K)^{2}(1-6K+\Delta \log (1-K)\big)\Big)=12K,$$
 or equivalently if and only if $a_{1}^{-}=0$ and either $a_{2}^{+}=a_{1}^{+}/2$ or $a_{2}^{-}=a_{1}^{+}/2$.

\end{theorem}
 
The condition $\Delta \log (1-K)=6K$ is equivalent to the flatness of the 
metric $(1-K)^{\frac{1}{3}}ds^2$, while the condition
 $$\Delta \log \big(\big(1-K)^{2}(1-6K+\Delta \log (1-K)\big)\big)=12K$$
 is equivalent to the flatness of the 
metric 
$$\big((1-K)^2(1-6K+\Delta \log (1-K))\big)^{\frac{1}{6}}ds^2.$$

The multiplication  on the Cayley numbers $\mathbb{O}$ yields a cross product on the purely imaginary Cayley numbers 
$\mathrm{Im}(\mathbb{O})=\mathbb{R}^7$ by
$$x\times y=\dfrac{1}{2}(x\cdot y-y\cdot x).$$
The scalar product on $\mathbb{R}^7$ is given by
$$\langle x,y\rangle=-\dfrac{1}{2}(x\cdot y+y\cdot x).$$
The almost complex structure on $\S^{6}$ is the endomorphism of its tangent bundle given by
$$J_xv=x \times v, {\ } x\in \S^{6}, v\in T_{x}\S^{6}.$$
$J$ is orthogonal and its covariant derivative is given by 
\begin{equation}\label{J}
 (\nabla_{X}J)Y=X \times Y +\langle X,JY\rangle x,
\end{equation} 
$X,Y$ being tangent vector fields. 

\begin{lemma}\label{B}
 Let $f:(M,ds^{2})\rightarrow \S^{6}$ be a pseudoholomorphic curve and $M$ be a Riemann surface with complex structure  $J$. For any vector fields $X,Y,Z$  we have:
$$B_{1}(JX,Y)=B_{1}(X,JY)=JB_{1}(X,Y),$$
\begin{eqnarray*}{}
 \nabla^{\bot}_{X}B_{1}(JY,Z)-J\nabla^{\bot}_{X}B_{1}(Y,Z)
=df(X)\cdot B_{1}(Y,Z)\\-df\circ J \circ A_{B_{1}(Y,Z)}X+df(A_{B_{1}(JY,Z)}X),
\end{eqnarray*}
where $A_{\xi}$ is the shape operator associated with a normal direction $\xi$ and $\nabla^{\bot}$ is the normal connection.
\end{lemma}

\begin{proof}
The lemma follows by differentiating twice $df\circ J= J\circ df$, using (\ref{J}), Gauss and Weingarten formulas and the fact that 
$df(X)\times df(Y)=-\langle X,JY\rangle f$.
\end{proof}

In particular, Lemma \ref{B} shows that pseudoholomorphic curves in $\S^6$ are superconformal 
surfaces.

\begin{lemma}\label{pseudoa}
 For every pseudoholomorphic curve $f:(M,ds^{2})\rightarrow \S^{6}$ we have $ 
a_{1}^{-}= 0$ and  $ a_{2}^{+}=a_{1}^{+}/2,$ or $ a_{2}^{-}=a_{1}^{+}/2$.
\end{lemma}
\begin{proof}
We choose an orthonormal frame along $f$ so that
$$B_{1}(e_{1},e_{1})=|B_{1}(e_{1},e_{1})|e_{3}, {\ } B_{1}(e_{1},e_{2})=|B_{1}(e_{1},e_{2})|e_{4},$$ 
$$e_{6}=df(e_{1})\cdot e_{3},{\ } e_{5}=Je_{6}.$$
From 
$$h_1^3=\kappa_{1}, \; h_2^3=0,\;  h_1^4=0,\; h_2^4=\kappa_{1}$$
we have 
$H_3=\kappa_{1}$, $H_4=i\kappa_{1}$ and thus $k_{1}^{+}=2\kappa_{1},
k_{1}^{-}= 0$.

We claim that $H_6=i(\kappa_{1}-H_5),$
or equivalently
$h^5_2=h^6_1$   and   $h^6_2=\kappa_{1}-h^5_1.$ Indded,  bearing
 in mind that the third fundamental form is given by
$$B_2(X,Y,Z)=\pi_2 \Big(\nabla^{\bot}_{X}B_1 (Y,Z)\Big),$$
where $\pi_2$ is the projection onto the second normal bundle, and using the second identity in Lemma \ref{B}  we have
$$h^5_2= \< J\nabla^{\bot}_{e_1}B_{1}(e_1,e_1) , e_5 \>=
\< \nabla^{\bot}_{e_1}B_{1}(e_1,e_1) , e_6 \> =h^6_1.$$
Similarly we have
$$h^6_2=\< B_{2}(e_1,e_1,e_2) , e_6 \>=\< B_{2}(e_1,Je_1,e_1) , e_6 \>=
\<\nabla^{\bot}_{e_1}B_{1}(Je_1,e_1) , e_6 \>.$$
Using the second identity in Lemma \ref{B}, and by the choice of the frame, we  
obtain
\begin{eqnarray*}
h^6_2\!\!\!&=&\!\!\!\<J\nabla^{\bot}_{e_1}B_{1}(e_1,e_1) , e_6 \>+
\< df(e_1) \cdot B_1(e_1,e_1),e_6 \>\\
\!\!\!&=&\!\!\!-\<J\nabla^{\bot}_{e_1}B_{1}(e_1,e_1) , Je_5 \>+
\kappa_1\< df(e_1) \cdot e_3,e_6 \>\\
\!\!\!&=&\!\!\!-\<\nabla^{\bot}_{e_1}B_{1}(e_1,e_1) , e_5 \>+
\kappa_1 \\
\!\!\!&=&\!\!\!
-h^5_1+\kappa_1.
\end{eqnarray*}
Thus,   $H_6=i(\kappa_{1}-H_5)$ , $k_{2}^{+}=\kappa_{1}, k_{2}^{-}= 2\overline{H}_5-\kappa_{1}$,
and the proof follows 
easily.
\end{proof}

\begin{lemma}\label{cI}
 Let $f:(M,ds^{2})\rightarrow \S^{6}$ be a  superconformal surface. 
If 
\begin{equation}\label{I}
\Delta \log (1-K)=6K-1, 
\end{equation}
then $f$ is superminimal with $a_2^+=a_1^+/2$. Conversely, if $f$ is superminimal and
 $a_2^+=a_1^+/2$, then (5.2) is satisfied.
\end{lemma}
\begin{proof}
We assume that (5.2) is satisfied and $f$ is not superminimal. Corollary \ref{superconformal} yields $b_{2}^{2}=1-K$ and 
 equations (\ref{supm}) become
$$\Delta \log a_{2}^{\pm} =3K\mp 2\dfrac{K_{2}^{\bot}}{b_{2}^{2}},$$
or equivalently on account of (5.2), $\Delta \log (1{\pm}F) =1\mp F$, where $F:=4K_{2}^{\bot}/b_{2}^{2}.$ From this we obtain
 $|\nabla F|^2=-(1-F^2)^2$,  which is a contradiction. Hence, $f$ is
superminimal and $a_2^+=a_1^+/2$ by Lemma \ref{pseudoa}.

Conversely, if we assume that $f$ is superminimal with $a_2^+=a_1^+/2$, then Corollary \ref{superconformal}
 implies (5.2).
\end{proof}

\begin{lemma}\label{III}
 Let $f:(M,ds^{2})\rightarrow \S^{6}$ be a  superconformal surface. 
If 
\begin{equation}\label{dI}
\Delta \log (1-K)=6K, 
\end{equation}
then $a_{2}^{+}=a_{2}^{-}=a_{1}^{+}/2$,  $a_{1}^{-}=0$ and $f$ lies in a totally geodesic $\S^5$ of $\S^6$.
Conversely, if 
$a_2^+=a_2^-=a_1^+/2$ and $a_{1}^{-}=0$, then (\ref{dI}) is satisfied and $f$ lies in a totally geodesic $\S^5$.
\end{lemma}
\begin{proof}
We assume that (\ref{dI}) is satisfied. Corollary \ref{superconformal} implies that $f$ cannot be superminimal. Then we have 
 $b_{2}^{2}=2(1-K)$ and 
 equations (\ref{supm}) become
$$\Delta \log a_{2}^{\pm} =3K\mp 4\dfrac{K_{2}^{\bot}}{b_{2}^{2}},$$
or equivalently, on account of (\ref{dI}), $\Delta \log (1{\pm}F) =\mp 2F$, where $F:=4K_{2}^{\bot}/b_{2}^{2}.$ We claim that 
$F$ is constant. Arguing indirectly, we assume that $\nabla F\neq 0$. Then from $\Delta F=-2F(1+F^2)$, 
$|\nabla F|^2=2F^{2}(1-F^2)$ and a well known argument\footnote{It is known (cf. \cite{EGT}) that
 if a two-dimensional Riemannian manifold
 $M$   allows a smooth  function $f:M\rightarrow 
\mathbb{R}$  such that $\Delta f=P(f)$ and $|\nabla
f|^{2}=Q(f)$ for smooth functions $P,Q:\mathbb{R}\rightarrow \mathbb{R},$
then on the set of points where the gradient $\nabla f\ $ doesn't vanish,  the Gaussian curvature $K$ satisfies 
\begin{equation*}
2KQ+(2P-Q^{\prime })(P-Q^{\prime })+Q(2P^{\prime }-Q^{\prime \prime })=0.
\end{equation*}},
 we  have that that $K=-8$, which contradicts (\ref{dI}). Hence, $F$ is constant and so $K_{2}^{\bot}=0$. 
This shows that $f$ lies in a totally geodesic $\S^5$ of $\S^6$ with $a_{2}^{+}=a_{2}^{-}=a_{1}^{+}/2$ and $a_{1}^{-}=0$.

Conversely, if $a_{2}^{+}=a_{2}^{-}=a_{1}^{+}/2$ and $a_{1}^{-}=0$, then $K_{2}^{\bot}=0$, $f$ lies in a totally geodesic $\S^5$ 
and Corollary \ref{superconformal}  implies (\ref{dI}).
\end{proof}

\begin{lemma}\label{II}
 Let $f:(M,ds^{2})\rightarrow \S^{6}$ be a non-superminimal  superconformal  surface. The condition
\begin{equation}\label{dII}
\Delta \log \Big(\big(1-K)^{2}(1-6K+\Delta \log (1-K)\big)\Big)=12K
\end{equation}
is satisfied if and only if either $a_2^+=a_1^+/2$ or $a_2^-=a_1^+/2$.
\end{lemma}
\begin{proof}
Assume that (\ref{dII}) is satisfied. From Corollary \ref{superconformal}, we have
\begin{equation*}
 b_2^2=(1-K)\big(2-6K+\Delta \log (1-K)\big).
\end{equation*} 
Then  (\ref{dII})  becomes
\begin{equation}\label{db}
 \Delta \log \big( \dfrac{b_2^2}{1-K}-1\big)=4-\dfrac{2b_2^2}{1-K}
\end{equation} 
and  equations (\ref{supm}) are written
$$\Delta \log u_{2}^{\pm} =2- u_{2}^{\pm},$$
where $u_{2}^{\pm}:=4(a_{2}^{\pm})^2/(1-K)$, or equivalently
\begin{equation}\label{u}
 \Delta u_{2}^{\pm}=\dfrac{|\nabla u_{2}^{\pm}|^2}{u_{2}^{\pm}}+u_{2}^{\pm}(2-u_{2}^{\pm}).
\end{equation} 
 Since $u_{2}^{+}+u_{2}^{-}=b_2^2/(1-K)$, (\ref{db}) is written as
$$\Delta \log (u_{2}^{+}+u_{2}^{-}-2)=4-(u_{2}^{+}+u_{2}^{-}),$$
or equivalently by virtue of (\ref{u}),
\begin{eqnarray*}{}
 u_{2}^{-}(u_{2}^{-}-2)|\nabla u_{2}^{+}|^2-2 u_{2}^{+}u_{2}^{-}\langle \nabla u_{2}^{+},\nabla u_{2}^{-} \rangle  
    +u_{2}^{+}(u_{2}^{+}-2)|\nabla u_{2}^{-}|^2\\
=2 u_{2}^{+}u_{2}^{-}(u_{2}^{+}+u_{2}^{-}-2)(u_{2}^{-}-2)(2-u_{2}^{+}).
\end{eqnarray*}
By the Cauchy-Schwartz inequality, we obtain
\begin{eqnarray*}{}
 u_{2}^{-}(u_{2}^{-}-2)|\nabla u_{2}^{+}|^2-2 u_{2}^{+}u_{2}^{-}|\nabla u_{2}^{+}||\nabla u_{2}^{-}|  
    +u_{2}^{+}(u_{2}^{+}-2)|\nabla u_{2}^{-}|^2 \\
\leq2 u_{2}^{+}u_{2}^{-}(u_{2}^{+}+u_{2}^{-}-2)(u_{2}^{-}-2)(2-u_{2}^{+})\\
\leq  u_{2}^{-}(u_{2}^{-}-2)|\nabla u_{2}^{+}|^2+2 u_{2}^{+}u_{2}^{-}|\nabla u_{2}^{+}||\nabla u_{2}^{-}|  
    +u_{2}^{+}(u_{2}^{+}-2)|\nabla u_{2}^{-}|^2.
\end{eqnarray*}
It easy to see that this double inequality holds only if $u_{2}^{+}=2$ or $u_{2}^{-}=2$. 
This yields $a_2^+=a_1^+/2$ or $a_2^-=a_1^-/2$.

Conversely, we assume that $a_2^+=a_1^+/2$ (the case  $a_2^-=a_1^+/2$ is similar). Corollary \ref{superconformal} 
 immediately implies  
 $$b_2^2=(1-K)\big(2-6K+\Delta \log (1-K)\big).$$
 From 
$$a^{\pm}_{2}:=\sqrt{\dfrac{1}{4}b_{2}^{2}\pm K_{2}^{\perp}}$$
and $a_2^+=a_1^+/2$, we find 
$$ K_{2}^{\perp}=\dfrac{1-K}{4}\Big(6K - \Delta \log (1-K) \Big)$$
and 
$$a_2^-= \dfrac{1}{\sqrt{2}}\sqrt{(1-K)\big(1-6K+\Delta \log (1-K) \big)}.$$
 Corollary \ref{superconformal}, implies that $\Delta \log (a_2^+a_2^-)=6K$, which yields (\ref{dII}).
\end{proof}

\begin{proof}[Proof of Theorem \ref{pseudo}.]
At first, we assume that $f$ is a pseudoholomorphic curve. According to Lemma \ref{pseudoa}, 
we have $a_{1}^{-}= 0,$ and either $ a_{2}^{+}=a_{1}^{+}/2$ or $ a_{2}^{-}=a_{1}^{+}/2$.

 If $f$ is superminimal, then $a_{2}^{-}= 0$, and  Lemma \ref{cI} yields
 $\Delta \log (1-K)=6K-1$.

 If $f$ lies in $\S^5$, then $K_2 ^{\bot}=0$, $a_{2}^{+}=a_{2}^{-}=a_{1}^{+}/2$ and Lemma \ref{III}  implies
 $\Delta \log (1-K)=6K$.

Assume that $f$ is neither superminimal nor lies in $\S^5$. From Corollary \ref{superconformal}, we have 
$$\|B_2\|^2=(1-K)\big(2-6K+\Delta \log (1-K)\big).$$
Assume further that $a_{2}^{+}=a_{1}^{+}/2$. From the proof of Lemma \ref{II}, we have $6K>\Delta \log (1-K)$ and  
$$4(a_{2}^{+}a_{2}^{-})^2=(1-K)^{2}(1-6K+\Delta \log (1-K)\big).$$
Corollary \ref{superconformal} implies that 
$$\Delta \log \Big(\big(1-K)^{2}(1-6K+\Delta \log (1-K)\big)\Big)=12K.$$

Now assume that $f$ is superconformal and satisfies one of (5.2), (\ref{dI}) or (\ref{dII}). In the case where (\ref{dII}) is fulfilled,  we further assume that 
$6K>\Delta \log (1-K)>6K-1$. Then the proof of the theorem follows from the preceding  lemmas and \cite{H}.
Indeed, our conditions imply that $(M,ds^2)$ satisfies the condition in \cite[Th. 6.1]{H}. Hence, there exists 
locally a pseudoholomorphic curve $g$ with induced metric $ds^2$.

If (5.2) is satisfied,
 or equivalently if $f$ is superminimal with $a_{2}^{+}=a_{1}^{+}/2$, then Lemma \ref{cI} implies that $g$ is also superminimal.
Since superminimal surfaces are rigid, we see that $f$ is $O(7)$-congruent to $g$.

If  (\ref{dI}) is satisfied,
 or equivalently if $a_{2}^{+}=a_{2}^{-}=a_{1}^{+}/2$ and $a_{1}^{-}=0$, then by Lemma \ref{III},  $f,g$ lie
 in a totally geodesic
 $\S^5$ and have the same $a$-invariants. From 
 Corollary \ref{superconformal}, we see that $f$ is $O(7)$-congruent to some  $g_{\theta}$.

Now assume that (\ref{dII}) is satisfied
 or equivalently  $a_{1}^{-}=0$ and either $a_{2}^{+}=a_{1}^{+}/2$ or $a_{2}^{-}=a_{1}^{+}/2$. Then Lemma \ref{II}  shows that
$f$ and $g$ have the same $a$-invariants, and the argument is the same as above.
\end{proof}

\section{The Ricci condition}
As an application of our main results, we provide another proof of the following result \cite{V2} that
 supports the Lawson's conjecture \cite{L}.
\begin{theorem}\label{v}
Lawson's conjecture is true for non-flat exceptional surfaces lying
 in odd-dimensional spheres.
\end{theorem}

We recall that  a two-dimensional Riemannian manifold $(M,ds^{2})$ with Gaussian curvature $%
K\leq 1$ satisfies the Ricci condition,
if and only  the metric $d\widehat{s}^{2}=\sqrt{1-K}ds^{2}$ is
flat away from points where $K=1,$ or equivalently  $\Delta \log (1-K)=4K.$

\begin{lemma}\label{v1}
 Let $f:(M,ds^{2})\rightarrow \S^{n}$ be a  non-flat exceptional surface
which satisfies the Ricci condition. Then the $a$-invariants of $f$  are given by 
\begin{equation*}
a_{r}^{+}=\left\{ 
\begin{array}{c}
c_{r}(1-K)^{\frac{r+1}{4}}\text{ \ if \ }r\text{ \ is odd,} \\ 
c_{r}(1-K)^{\frac{r+2}{4}}\text{ \ if \ }r \equiv 2 {\ }\mathrm{mod}{\ }4,\\ 
c_{r}(1-K)^{\frac{r}{4}}\text{ \ if \ }r \equiv 0 {\ }\mathrm{mod}{\ }4,\\
\end{array}%
\right.
\end{equation*}%
\begin{equation*}
 a^{-}_{r}=a^{+}_{r}\sqrt{\frac{1-\rho _{r}}{1+\rho _{r}}},
\end{equation*}
where $c_{r}=2^{\frac{2-r}{2}}\beta_{r}$, for $0\leq r\leq m-1,$  $m=[(n-1)/2],$ $c_{m}=2^{\frac{1-m}{2}}\beta_{m},$
 $\beta_{1}=\sqrt{2}, \beta_{r} = \rho _{r-1}\beta_{r-1}$ 
 and $\rho _{r}=1$ if $r$ is even.
\end{lemma}
\begin{proof}
 The lemma follows easily by induction using the Ricci condition and Theorem 2.
\end{proof}

\begin{proof}[Proof  of Theorem \ref{v}.]
We claim that 
$n\equiv 3{\ } {\mathrm {mod}}{\ } 4.$ Arguing indirectly, we suppose that $n=4m+1.$ Then Lemma \ref{v1}
 yields $a^+_{2m}=c_{2m}(1-K)^{\frac{m}{2}}$ if $m$ is even and $a^+_{2m}=c_{2m}(1-K)^{\frac{m+1}{2}}$ if $m$ is odd. 
Moreover, viewing $%
f$ as a minimal surface in $\S^{4m+2},$ we obviously have $K_{2m}^{\perp
}=0.$ Then from Theorem 2, we obtain $\Delta \log
a^+_{2m}=(2m+1)K.$ Then the Ricci condition yields $K=0,$ which is a
contradiction.

Hence $n=4m+3.$ According to Lemma \ref{v1}, $\Phi _{r}=0$ if $r$ is
even. Let $r_{0}=\min \left\{ r:1\leq r\leq 2m+1\text{ with }\Phi _{r}\neq
0\right\}.$ Obviously $r_{0}$ is odd. Let $z$ be a local complex coordinate
such that $ds^{2}=F|dz|^{2}.$  We pick a branch $g$ of $%
f_{r_{0}}^{\frac{2}{r_{0}+1}}$, where $f_{r}=\langle
B_{r}^{(r+1,0)},B_{r}^{(r+1,0)}\rangle,$ and define the quadratic form $\Phi =gdz^{4}.$
It is obvious that $\Phi $ is well defined and holomorphic. For any odd $%
r\geq r_{0},$ we write $\Phi _{r}=|f_{r}|e^{i\tau _{r}}dz^{2r+2}.$
Appealing to Lemma \ref{v1} and (\ref{hl}), we obtain
\begin{equation*}
\Phi _{r}=\gamma _{r}e^{i\big(%
\tau _{r}-\frac{r+1}{r_{0}+1}\tau _{r_{0}}\big)}\Phi ^{\frac{r+1}{2}},
\end{equation*}%
where $\gamma _{r}$ is a positive number. 
From the holomorphicity of $\Phi _{r}$ and $\Phi,$ we deduce that $\tau
_{r}-\frac{r+1}{r_{0}+1}\tau _{r_{0}}$ is constant. Moreover, we 
see that $|g|^{2}=c_{0}F^{4}(1-K)^{2},$ where $c_{0}$ is a positive
constant. Using the holomorphicity of $g$ and arguing as in \cite[Theorem 8]%
{L0}, we infer that\ there exists locally a minimal surface $\widetilde{f}$
in $\S^{3}$ with Hopf differential $\widetilde{\Phi }=c\Phi,$ where $c$ is a
complex number. Therefore $\Phi _{r}=\delta_{r}\widetilde{\Phi }^{\frac{r+1}{2}}$
for any odd $r\geq r_{0},$ where $\delta_{r}$ is a complex number and $\Phi
_{r}=0 $ otherwise.\ By \cite[Prop. 2]{V2}, the Hopf differentials
of $f$\ coincide with those of minimal surfaces which decompose as a direct
sum of the associated minimal surfaces in $\S^{3}.$ According to the main
result in \cite{V},  $f$ splits as a direct sum of the associated
minimal surfaces of $\widetilde{f}.$
\end{proof}

\section{Global formulas}
In this section, we give some topological restrictions for exceptional surfaces. 
The zero set  of an AVT function $a$ on a connected compact oriented surface $M$ is either isolated or the whole of $M$, and outside 
its zeros, the function is smooth. 
If $a$ is a non-zero AVT function, i.e., locally $a=|t_{0}|a_{1}$, with $t_0$ holomorphic, the order $k\geq 1$ of any $p \in M$
 with $a(p)=0$ 
is the order of $t_0$ at $p$. Let $N(a)$ be the sum of all orders for all zeros of $a$. Then $\Delta \log a$ is bounded on
 $M\smallsetminus \{ a=0  \}$ and its integral is computed in the following lemma which was proved in \cite{EGT,ET}.

\begin{lemma}
Let $(M,ds^{2})$ be a compact oriented two-dimensional Riemannian manifold
with area element $dA.$ If $a$ is an AVT function on $M,$ then 
\begin{equation*}
\int_{M}\Delta \log adA=-2\pi N(a).
\end{equation*}
\end{lemma}

For exceptional surfaces it has been proved in \cite[Prop. 4]{V2} that all higher normal bundles can 
be smoothly extended over the whole surface. Then the following follows from Lemma 10,  Theorem 2 and Proposition \ref{Asp}.
\bigskip

\begin{corollary}
Let $f:(M,ds^{2})\rightarrow \S^{n}$ be an  exceptional surface. The Euler number $\chi (N^{r}f)$ 
of the $r$-th
normal bundle and the Euler-Poincar\'{e} characteristic $\chi (M)$ of $M$ satisfy the following:

(i) If $\Phi_{r}\neq 0$ for some $1\leq r< m$, where $m=[(n-1)/2]$, then 
$$  \chi (N^{r}f)=0  {\ }  \mathrm{and}   {\ }(r+1)\chi (M)=-N(a^{+}_r)=-N(a^{-}_r).$$

(ii) If $\Phi_{r} = 0,$ for some $1\leq r \leq m$, then 
$$(r+1)\chi (M)- \chi (N^{r}f)=-N(a^{+}_r).$$

(iii)  If $\Phi_{m}\neq 0,$ then 
$$(m+1)\chi (M)\mp \chi (N^{m}f)=-N(a^{\pm}_m).$$
\end{corollary}

\bigskip

\end{document}